\documentclass[draft]{amsart}

\usepackage{amscd}

\newtheorem{thm}{Theorem}
\newtheorem{prop}[thm]{Proposition}
\newtheorem{lem}[thm]{Lemma}
\newtheorem{cor}[thm]{Corollary}

\theoremstyle{remark}
\newtheorem{rmk}[thm]{Remark}

\theoremstyle{definition}
\newtheorem{dfn}[thm]{Definition}

\numberwithin{thm}{section}

\begin{document}


\title{Attractors on $\mathbf{P}^k$}

\author{Feng Rong}

\address{Department of Mathematics, University of Michigan, Ann Arbor, MI 48109, USA}
\email{frong@umich.edu}

\begin{abstract}
We show that special perturbations of a particular holomorphic map on $\mathbf{P}^k$ give us examples of maps that possess chaotic nonalgebraic attractors.  Furthermore, we study the dynamics of the maps on the attractors. In particular, we construct invariant hyperbolic measures supported on the attractors with nice dynamical properties.
\end{abstract}

\maketitle

\section{Introduction}

Attractors in holomorphic dynamics have been studied in \cite{FW:Attractor}, \cite{FS:Examples} and \cite{JW:Nonalgebraic}.
In \cite{FW:Attractor}, Forn$\ae$ss and Weickert gave some basic properties of attractors. In \cite{JW:Nonalgebraic}, Jonsson and Weickert gave
a nonalgebraic attractor on $\mathbf{P}^2$. In this paper, we extend their results to arbitrary dimensions and explore further dynamical properties
of the map on the attractor. More precisely, we have the following.

\begin{thm}\label{T:Main1}
Let $f_\lambda$ be a holomorphic map on $\mathbf{P}^k$, $k\ge2$, of the form
$$f_\lambda[z:w_1:\cdots:w_{k-1}:t] = [(z-2w_1)^2:\cdots:(z-2w_{k-1})^2:z^2:t^2+\lambda z^2].$$
The map $f_\lambda$ has a nonalgebraic attractor $K_\lambda$ for $\lambda\in \mathbf{C}$ with $|\lambda|\neq 0$ sufficiently small. Moreover, 
the map $f_\lambda$ is chaotic on $K_\lambda$.
\end{thm}

\begin{thm}\label{T:Main2}
Let $f_\lambda$ and $K_\lambda$ be as in Theorem \ref{T:Main1}. There exists a probability measure $\mu_\lambda$ with 
support equal to $K_\lambda$ and invariant under $f_\lambda$, i.e. ${f_\lambda}_\star \mu_\lambda = \mu_\lambda$, such that\\
(1) The measure $\mu_\lambda$ is mixing;\\
(2) The measure $\mu_\lambda$ is the unique measure of maximal entropy $(=(k-1)\log2)$ for $f_\lambda|_{K_\lambda}$;\\
(3) The measure $\mu_\lambda$ describes the distribution of periodic points for $f_\lambda|_{K_\lambda}$;\\
(4) The smallest nonnegative Lyapunov exponent of $f_\lambda|_{K_\lambda}$ with respect to $\mu_\lambda$ is greater or equal to 
$\frac{1}{2}\log2$ at $\mu_\lambda-$almost every point.
\end{thm}

The structure of the paper is as follows. In \S \ref{S:Prelim} we recall some facts from ergodic theory and holomorphic dynamics and give
a precise definition of an attractor. As a necessary preparation, we study the so-called history space in \S \ref{S:History}, which is also of its own interest. 
Theorem \ref{T:Main1} is proved in \S \ref{S:Attractor} and  Theorem \ref{T:Main1} in \S \ref{S:Measure}. 

\section{Preliminaries}\label{S:Prelim}

The material in this section is fundamental and will be needed in consequent sections. In \S \ref{SS:Erg} we recall some notions and results 
from ergodic theory. In \S \ref{SS:Hol} we give a brief account on holomorphic dynamical systems. In \S \ref{SS:Attr} we give a precise definition
of an attractor.

\subsection{Ergodic Theory}\label{SS:Erg}

Let $M$ be a compact metric space with a metric $d$ and $f:M\rightarrow M$ be a continuous map.

\begin{dfn}
Let $U,V$ be two nonempty open sets in $M$.  Then $f$ is $topologically$ $transitive$ on $M$, if there exists $n\in\mathbf{N}$
such that $f^n(U)\cap V \neq \emptyset$. And $f$ is $topologically$ $mixing$ on $M$, if there exists $N > 0$ such that 
$f^n(U)\cap V \neq \emptyset$ for all $n\ge N$.
\end{dfn}

An important concept in ergodic theory is $entropy$ (see \cite{KH:Intro},\cite{W:Ergodic} for details). For any $n\in \mathbf{N}$,
we can define a new metric on $M$ by 
$$d_n^f(x,y)=\max_{0\le i\le n-1}d(f^i(x),f^i(y)).$$ 
A subset $N$ of $M$ is said to $(n,\epsilon)-span$ $M$ with respect to $f$ if for any $x\in M$ there exists $y\in N$ with 
$d_n^f(x,y)\le \epsilon$, $\epsilon>0$.

\begin{dfn}
The $topological\ entropy$ of $f$ is
$$h(f)=\lim\limits_{\epsilon\rightarrow 0}\limsup\limits_{n\rightarrow\infty}\frac{1}{n}\ \log\ r_n^f(\epsilon,M),$$
where $r_n^f(\epsilon,M)$ is the smallest cardinality of any $(n,\epsilon)$-spanning set for $M$ with respect to $f$.
\end{dfn}

Let $\mu$ be a probability measure on $M$ such that $f$ is measure-preserving, i.e. for any measurable set $A\subset M$, $\mu(f^{-1}(A))=\mu(A)$.

\begin{dfn}
Let $A$ and $B$ be any two measurable subsets of $M$.\\ 
$(i)$ The map $f$ is $ergodic$ if
$$\lim\limits_{n\rightarrow\infty}\frac{1}{n}\sum\limits_{i=0}^{n-1}\mu(f^{-i}(A)\cap B)=\mu(A)\mu(B);$$
$(ii)$ The map $f$ is $weak-mixing$ if
$$\lim\limits_{n\rightarrow\infty}\frac{1}{n}\sum\limits_{i=0}^{n-1}|\mu(f^{-i}(A)\cap B)-\mu(A)\mu(B)|=0;$$
$(iii)$ The map $f$ is $mixing$ if
$$\lim\limits_{n\rightarrow\infty}\mu(f^{-n}(A)\cap B)=\mu(A)\mu(B).$$
\end{dfn}

Let $P(M,f)$ denote the set of probability measures on $M$ such that $f$ is measure-preserving. For each $\nu\in P(M,f)$
we can define the $measure-theoretic$ $entropy$ $h_\nu(f)$. The following theorem by M. Brin and A. Katok will serve as the definition
of the measure-theoretic entropy for our purpose, which will also be needed later.

\begin{thm}\cite{BK:Entropy}\label{T:BK}
For $\nu$-almost every $x\in M$\\
(a) $h_\nu(f,x)\ :=\ \lim\limits_{\epsilon\rightarrow 0}\liminf\limits_{n\rightarrow\infty}-\frac{1}{n} \log \nu(B_n^f(x,\epsilon))
\ =\ \lim\limits_{\epsilon\rightarrow 0}\limsup\limits_{n\rightarrow\infty}-\frac{1}{n} \log \nu(B_n^f(x,\epsilon))$,\\
where $B_n^f(x,\epsilon)$ is the $d_n$-ball about $x$ of radius $\epsilon$;\\
(b) $h_\nu(f,x)$ is f-invariant;\\
(c) $h_\nu(f)=\int\limits_M h_\nu(f,x)\ d\nu$.
\end{thm}

The Variational Principle says that $h(f)=sup\{h_\nu(f) | \nu\in P(M,f)\}$ (see \cite{W:Ergodic}, p.187). A measure that attains the
supremum, if it exists, will be called a measure of maximal entropy.

For $n\in \mathbf{N}$, let $E_n(f)=\{x\in M | f^n(x)=x\}$, i.e. the set of periodic points whose periods are
factors of $n$, and let $N_n(f)$ denote the cardinality of $E_n(f)$. A measure $\nu\in P(M,f)$ is said to describe
the distribution of periodic points for $f$ if
$$\frac{1}{N_n(f)}\sum\limits_{x\in E_n(f)}\ \delta_x\ \longrightarrow \nu\ \ weakly.$$

\subsection{Holomorphic Dynamical Systems}\label{SS:Hol}

Let $f=[f_0:\cdots:f_k]:\mathbf{P}^k\rightarrow \mathbf{P}^k$ be a holomorphic map on $\mathbf{P}^k$ of (algebraic) degree $d\ge 2$,
i.e. each $f_j$ is a homogeneous polynomial of degree $d$. Let $F=(f_0,\cdots,f_k)$ be its lifting to $\mathbf{C}^{k+1}$. 
The dynamical Green function on $\mathbf{C}^{k+1}$ is defined as 
$$G(z)=\lim\limits_{n\rightarrow\infty}\frac{1}{d^n}\log|F^n(z)|,$$
where $|\cdot|$ denotes the Euclidean norm on $\mathbf{C}^{k+1}$. It is continuous and plurisubharmonic on $\mathbf{C}^{k+1}\backslash \{0\}$
and satisfies $G\circ F=dG$ (\cite{FS:Potential}). The Green current $T$ associated with $f$ is then defined by 
$\pi^{*}T=dd^cG$, where $\pi$ is the canonical projection from $\mathbf{C}^{k+1}\backslash \{0\}$ to $\mathbf{P}^k$. It is a positive closed 
$(1,1)$ current on $\mathbf{P}^k$ with $\| T \|=1$. For more information on currents, see \cite{FS:Potential}, \cite{deR:Current}.

The significance of the Green current $T$ lies in the fact that the support of $T$ is equal to the Julia set $J_1$ of $f$, i.e. the complement
of the Fatou set $\Omega$. The Fatou set $\Omega$ consists of points each of which has a neighborhood where the family 
$\{f^j\}_{j=1}^{\infty}$ is a normal family. It is shown in \cite{U:Fatou} that $\pi^{-1}(\Omega)$ is equal to the set where the Green
function $G$ is pluriharmonic. 

For each $l\le k$, we can define $T^l$ as the $l-$th wedge product of $T$.  Locally on an open set $U\subset \mathbf{P}^k$, $T^l:=(dd^c(G\circ s))^l$, 
where $s$ is a holomorphic section of $\pi$ on $U$. The definition makes sense since $G$ is continuous and plurisubharmonic on
$\mathbf{C}^{k+1}\backslash \{0\}$ (\cite{BT:Capacity}). We can then define the intermediate Julia sets $J_l$ to be the support of $T^l$. 
It is obvious that each $J_l$ is totally invariant and $J_1\supset J_2\supset \cdots \supset J_k$. We also have 
$f^{*}(T^l)=d^lT^l$, $f_{*}(T^l)=T^l$ and $\|T^l \|=1$.

The current $\mu:=T^k$ is of particular importance. It is an $f-$invariant probability measure on $\mathbf{P}^k$. It is mixing and is the
unique measure of maximal entropy ($=k\log d$), see \cite{FS:Potential}, \cite{BrD:Deux}. It describes the distribution of periodic points for $f$ 
(\cite{BrD:Exposants}), i.e.
$$\frac{1}{d^{kn}}\sum\limits_{z\in E_n} \delta_z\longrightarrow \mu,\ \ \ weakly.$$
In particular, repelling periodic points are dense in $J_k=supp(\mu)$. It also describes the distribution of preimages of points 
outside of a pluripolar set (\cite{BrD:Deux},\cite{FS:Potential}), i.e.
$$\frac{1}{d^{kn}}\sum\limits_{f^n(w)=z} \delta_w\longrightarrow \mu,\ \ \ weakly,\ \ \ z\in \mathbf{P}^k\backslash E,\ \ \ E\ pluripolar.$$ 

\subsection{Attractors}\label{SS:Attr}

Let $M$ be a compact metric space with a metric $d$ and $f: M\rightarrow M$ be a continuous map.

\begin{dfn}(\cite{R:Dynamics})
A sequence $(x_i)_{0\le i\le n}$ of points in $M$ is an $\epsilon-pseudoorbit$ if $d(f(x_i),x_{i+1}) < \epsilon$ for $i = 0,\cdots,n-1$, $\epsilon>0$.
For $x,y\in M$, we write $x\succ y$ if for every $\epsilon>0$ there exists an $\epsilon-pseudoorbit$ $(x_i)_{0\le i\le n}$ with $x_0 = x$ and $x_n = y$. 
The preorder $\succ$ is transitive. We write $x\sim y$ if $x\succ y$ and $y\succ x$. This defines an equivalence relation on $M$, the equivalence 
classes of which are closed subsets of $M$. The ordering $\succ$ induces an ordering between equivalence classes. An $attractor$ for $f$ is 
a minimal equivalence class.
\end{dfn}

There is also the notion of an $attracting \ set$.

\begin{dfn}
A nonempty compact subset $K\subset M$ is an $attracting$ $set$ if $K$ has a neighborhood $U$ such that $f(U)\subset\subset U$ and 
$K = \bigcap_{n\ge 0} f^n(U)$.
\end{dfn}

\begin{lem}
If $f$ is topologically transitive on $M$, then there is a dense subset $R$ of $M$ such that for each $x\in R$ its orbit $\{f^n(x):n\ge 0\}$ is dense in $M$.
\end{lem}
\begin{proof}
Let $\{U_i\}_{i\in \mathbf{N}}$ be a countable basis for $M$. Then every point $x$ in the set $R := \cap_{i=1}^\infty\cup_{n=0}^\infty f^{-n}(U_i)$ 
has a dense orbit. Suppose not, then there exists $x\in R$ such that $\cup_{n=0}^{\infty}f^n(x)$ misses an open set $U$. Let $U_1\subset U$ 
be a basis open set, we have $x\notin \cup_{n=0}^{\infty}f^{-n}(U_1)$, a contradiction. Note that $\cup_{n=0}^\infty f^{-n}(U_i)$ is open and dense 
in $M$ for every $i\in \mathbf{N}$ by the continuity and topological transitivity of $f$. Hence $R$ is dense in $M$ by Baire's theorem.
\end{proof}

\begin{rmk}\label{R:Attractor}
It follows easily from the previous lemma that if $K$ is an attracting set for $f$ and $f$ is topologically transitive on $K$, then $K$ is an attractor for $f$.
 \end{rmk}

\begin{rmk}
Milnor gave another definition of an attractor in \cite{M:Attractor} from a measure-theoretic point of view. 
\end{rmk}

\begin{dfn}(\cite{D:Chaotic})
We say that $f$ has sensitive dependence on initial conditions if there exists $\delta>0$ such that the following holds.
For any $x\in M$ and any neighborhood $U$ of $x$, there exist $y\in U$ and $n\ge 0$ such that $d(f^n(x),f^n(y))>\delta$.
\end{dfn}

\begin{dfn}(\cite{D:Chaotic})
We say that $f$ is $chaotic$ on $M$ if\\
$(i)$ The map $f$ is topologically transitive on $M$;\\
$(ii)$ periodic orbits for $f$ are dense in $M$;\\
$(iii)$ The map $f$ has sensitive dependence on initial conditions.
\end{dfn}

\begin{rmk}
See \cite{BBCDS} for an interesting discussion of this definition.
\end{rmk}

\begin{prop}
If $f$ is topologically mixing on $M$, then $f$ has sensitive dependence on initial conditions.
\end{prop}
\begin{proof}
Fix $\delta>0$ small enough so that we can find two open sets $V_1$ and $V_2$ satisfying $\tilde{V_1}\cap \tilde{V_2} = \emptyset$,
where $\tilde{V_i}=\{ x\in M | d(x,V_i)\le \delta \}$, $i=1,2$. Since $f$ is topologically mixing on $M$, for any $x\in M$ and any neighborhood
$U$ of $x$, there exists $N>0$ such that $f^n(U)\cap V_i\neq \emptyset$ for all $n\ge N$. If for some $n\ge N$ 
$f^n(x)\in M\backslash \tilde{V_1}$, then for any $x_1\in f^n(U)\cap V_1$ and $y_1\in U$ with $f^n(y_1)=x_1$ we have
$d(f^n(x), f^n(y_1)) > \delta$. Otherwise the set $W=\cup_{n\ge N} f^n(x)$ satisfies $W\subset \tilde{V_1}$. But then
$W\subset M\backslash \tilde{V_2}$ since $\tilde{V_1}\cap \tilde{V_2} = \emptyset$. Hence for any $n\ge N$, any point 
$x_2\in f^n(U)\cap V_2$ and $y_2\in U$ with $f^n(y_2)=x_2$ we have $d(f^n(x), f^n(y_2)) > \delta$. This completes the proof.
\end{proof}

\begin{rmk}\label{R:Mixing}
Combining the above proposition and Remark \ref{R:Attractor}, we see that if $f$ is topologically mixing on
an attracting set $K$ then $K$ is a chaotic attractor provided that periodic points for $f$ are dense in $K$.
\end{rmk}

\section{History Space}\label{S:History}

The notion of history space is very important in the dynamical study of endomorphisms. It is of its own interest besides serving as
a necessary preparation for our study of attractors. Therefore we state our results in a more general setting. 

Let $M$ be a compact metric space and $f: M\rightarrow M$ be a continuous map. 

\begin{dfn}
The $history\ space$ $\hat{M}$ associated with the map $f$ on $M$ is the set
$$\hat{M} := \left\{\{x_{-i}\}\in \prod_{i=0}^{\infty} M | f(x_{-(i+1)})=x_{-i}, \forall\ i\ge 0\right\}.$$
\end{dfn}

There is a collection of canonical projections $\{\pi_j:\hat{M}\rightarrow M\}_{j\ge 0}$, where each $\pi_j$ projects the history 
$\{x_{-i}\}$ to its member $x_{-j}$. There is also a natural lifted map $\hat{f}:\hat{M}\rightarrow \hat{M}$ of $f$, which maps the history 
$\{x_{-i}\}$ to the history $\{y_{-i}\}$ with $y_{-i}=f(x_{-i})$ for $i\ge 0$. The following relationship is clear:
$$f\circ \pi_j = \pi_j\circ \hat{f}\ \ \ \ \forall\ j\ge 0.$$

The standard topology for the history space $\hat{M}$ is obtained by considering it as a subset
of the compact space $\prod_{i=0}^{\infty} M$ with the usual product topology. As a closed subspace
of a compact space, $\hat{M}$ is compact in this topology. Note that each canonical projection
$\pi_j$ is continuous. Let $d$ be a metric giving the topology of $M$ and bounded by 1, then 
a metric $\hat{d}$ giving the topology of $\hat{M}$ and bounded can be defined as follows.
$$\hat{d}(\hat{x},\hat{y})=\sum\limits_{i=0}^{\infty}2^{-i}\ d(x_{-i},y_{-i}),\ \ \hat{x}=\{x_{-i}\},\ \hat{y}=\{y_{-i}\}.$$

\begin{lem}\label{L:Basis}
Let $B$ be a basis of open sets in $M$. Then $\hat{B}=\{\pi_j^{-1}(U) | U\in B,\ j\ge 0\}$ is a basis of open sets in $\hat{M}$.
\end{lem}
\begin{proof}
As a subspace of $\prod_{i=0}^{\infty} M$, any open set in $\hat{M}$ is of the form
$\hat{M}\cap V$, where $V$ is an open set in $\prod_{i=0}^{\infty} M$. By the definition of the
product topology on $\prod_{i=0}^{\infty} M$, the set
$$A=\{\prod_{i=0}^{\infty} V_i | V_i\ open\ in\ M,\ V_i=M\ for\ all\ but\ finitely\ many\ i's\}$$
is a basis of open sets in $\prod_{i=0}^{\infty} M$. Therefore it suffices to show that for any $V\in A$, the set $\hat{M}\cap V$ has the form 
$\pi_j^{-1}(U)$ for some $U$ open in $M$ and $j\ge 0$. Assume $V=\prod_{i=0}^{\infty} V_i$ with $V_i\neq M$ for $i\in \{i_1,\cdots,i_n\}$ and
$i_1<\cdots < i_n$. Then by the definition of $\hat{M}$, we have $\hat{M}\cap V=\pi_{i_n}^{-1}(U)$, where $U=\cap_{j=1}^{n} f^{i_j-i_n}(V_{i_j})$ 
is an open set in $M$. This completes the proof.
\end{proof}

\begin{lem}
The map $\hat{f}:\hat{M}\rightarrow \hat{M}$ is a homeomorphism.
\end{lem}
\begin{proof}
Let $\hat{U}$ be any open set in $\hat{M}$. Without loss of generality, we can write $\hat{U}=\pi_j^{-1}(U)$ with $U$ open in $M$ and $j\ge 0$.
By the continuity of both $f$ and $\pi_j$, $\hat{f}^{-1}(\pi_{j}^{-1}(U))=\pi_{j}^{-1}(f^{-1}(U))$ is open. Hence $\hat{f}$ is continuous. Since 
$\hat{f}(\pi_{j}^{-1}(U))=\pi_{j+1}^{-1}(U)$ is open, the map $\hat{f}^{-1}$ is also continuous. Therefore $f$ is a homeomorphism as it is also 
one-to-one and onto.
\end{proof}

We now begin the study of the dynamical properties of $\hat{f}$ on $\hat{M}$. We will see that the dynamics of
$\hat{f}$ on $\hat{M}$ is closely related to that of $f$ on $M$.

\begin{thm}\label{T:Transitive}
The map $\hat{f}$ is topologically transitive (resp. mixing) on $\hat{M}$ if and only if the map $f$ is topologically transitive 
(resp. mixing) on $M$.
\end{thm}
\begin{proof}
Let $U$ and $V$ be two open sets in $M$. Then $\pi_0^{-1}(U)$ and $\pi_0^{-1}(V)$ are two open sets
in $\hat{M}$. If $\hat{f}$ is topologically transitive on $\hat{M}$, then there exists $n>0$
such that $\hat{f}^n(\pi_0^{-1}(U))\cap \pi_0^{-1}(V)\neq \emptyset$. Hence $f^n(U)\cap V\neq \emptyset$, i.e.
$f$ is topologically transitive. If $\hat{f}$ is topologically mixing on $\hat{M}$, then there exists 
$N>0$, such that for any $n\ge N$ $\hat{f}^n(\pi_0^{-1}(U))\cap \pi_0^{-1}(V)\neq \emptyset$. Hence 
$f^n(U)\cap V\neq \emptyset$, i.e. $f$ is topologically mixing.

Conversely, let $\hat{U}$ and $\hat{V}$ be two open sets in $\hat{M}$. Without loss of generality, we can write $\hat{U}=\pi_i^{-1}(U)$ and
$\hat{V}=\pi_j^{-1}(V)$ with $U$ and $V$ open in $M$ and $i,j\ge 0$. 
If $f$ is topologically transitive on $M$, then there exists $n>0$, such that $f^n(f^{i}(U))\cap f^{j}(V)\neq \emptyset$. Hence
$\hat{f}^n(\pi_i^{-1}(U))\cap \pi_j^{-1}(V)\neq \emptyset$, i.e. $\hat{f}$ is topologically transitive.
If $f$ is topologically mixing on $M$, then there exists $N>0$ such that for any $n\ge N$,
$f^n(f^{i}(U))\cap f^{j}(V)\neq \emptyset$. Hence $\hat{f}^n(\pi_i^{-1}(U))\cap \pi_j^{-1}(V)\neq \emptyset$,
i.e. $\hat{f}$ is topologically mixing.
\end{proof}

Let $\mu$ be a probability measure on $M$ such that $f$ is measure-preserving. We can define an induced
measure $\hat{\mu}$ on $\hat{M}$ as follows.

\begin{dfn}\label{D:Measure}
The $history\ measure$ $\hat{\mu}$ on $\hat{M}$ induced by the measure $\mu$ on $M$ is the
measure $\hat{\mu}(\hat{A})=\lim\limits_{j\rightarrow\infty}\mu(\pi_j(\hat{A}))$,
for any subset $\hat{A}\subset \hat{M}$ such that $\pi_j(\hat{A})$ is $\mu-$measurable for each $j\ge 0$.
\end{dfn}

\begin{rmk}\label{R:Limit}
Since $f$ is measure-preserving and $\pi_{j+1}(\hat{A})\subset f^{-1}(\pi_j(\hat{A}))$ for any $j\ge 0$, $\{\mu(\pi_j(\hat{A}))\}_{j=0}^{\infty}$
is a non-increasing sequence of positive numbers. Therefore, the limit $\lim\limits_{j\rightarrow\infty}\mu(\pi_j(\hat{A}))$ always exists.
\end{rmk}

\begin{rmk}\label{R:History}
For a full history $\hat{A}=\pi_0^{-1}(A),\ A\subset M$ measurable, we have $\hat{\mu}(\hat{A})=\mu(A)$. This is
because $\mu(\pi_j(\hat{A}))=\mu(A)$ for any $j>0$, as $f$ is measure-preserving. In particular, 
$\hat{\mu}(\hat{M})=1$, proving that $\hat{\mu}$ is a probability measure on $\hat{M}$.
\end{rmk}

\begin{lem}
For any subset $\hat{A}\subset \hat{M}$, we have 
$\hat{\mu}(\hat{f}(\hat{A}))=\hat{\mu}(\hat{f}^{-1}(\hat{A}))=\hat{\mu}(\hat{A})$. 
In particular, the map $\hat{f}$ is measure-preserving.
\end{lem}
\begin{proof}
It follows from the definition.
\end{proof}

Due to the nature of the definition of the history measure, one might expect it to have similar properties as that of the measure on the 
base space. Indeed, this turns out to be true. We are going to show some of them in the sequel. 

\begin{thm}\label{T:Ergodic}
The map $\hat{f}$ is ergodic (resp. weak-mixing, mixing) relative to $\hat{\mu}$ on $\hat{M}$ if and only if 
the map $f$ is ergodic (resp. weak-mixing, mixing) relative to $\mu$ on $M$.
\end{thm}
\begin{proof}
Let $A$, $B$ be any two measurable subsets of $M$ and $\hat{A}=\pi_0^{-1}(A)$, $\hat{B}=\pi_0^{-1}(B)$
be their full histories. By Remark \ref{R:History}, we have $\hat{\mu}(\hat{A})=\mu(A)$, $\hat{\mu}(\hat{B})=\mu(B)$ and
$\hat{\mu}(\hat{f}^{-1}(\hat{A})\cap \hat{B})=\mu(f^{-1}(A)\cap B)$. Therefore, if $\hat{f}$ is ergodic 
(resp. weak-mixing, mixing) then $f$ is ergodic (resp. weak-mixing, mixing).

Conversely, let $\hat{A}$, $\hat{B}$ be any two measurable subsets of $\hat{M}$. By Definition \ref{D:Measure} and
Remark \ref{R:Limit}, for any $\epsilon>0$, there exists an integer $N\ge 0$ such that for any $l\ge N$ we have
$$\mu(\pi_l(\hat{A}))-\hat{\mu}(\hat{A})<\epsilon,\ \  \mu(\pi_l(\hat{B}))-\hat{\mu}(\hat{B})<\epsilon$$
and
$$|\mu(f^{-l+N}(\pi_N(\hat{A}))\cap\pi_N(\hat{B}))-\hat{\mu}(\hat{f}^{-l+N}(\hat{A})\cap\hat{B})|<\epsilon.$$

For any $j\ge 0$ we have
\begin{eqnarray*}
& &\hat{\mu}(\hat{f}^{-j}(\hat{A})\cap \hat{B})-\hat{\mu}(\hat{A})\hat{\mu}(\hat{B})\\
&=&\hat{\mu}(\hat{f}^{-j}(\hat{A})\cap \hat{B})-\mu(f^{-j}(\pi_N(\hat{A}))\cap \pi_N(\hat{B}))\\
& &+\mu(f^{-j}(\pi_N(\hat{A}))\cap \pi_N(\hat{B}))-\mu(\pi_N(\hat{A}))\mu(\pi_N(\hat{B}))\\
& &+\mu(\pi_N(\hat{A}))\mu(\pi_N(\hat{B}))-\hat{\mu}(\hat{A})\hat{\mu}(\hat{B}).
\end{eqnarray*}
Hence
\begin{eqnarray*}
& &|[\hat{\mu}(\hat{f}^{-j}(\hat{A})\cap \hat{B})-\hat{\mu}(\hat{A})\hat{\mu}(\hat{B})]\\
& &-[\mu(f^{-j}(\pi_N(\hat{A}))\cap \pi_N(\hat{B}))-\mu(\pi_N(\hat{A}))\mu(\pi_N(\hat{B}))]|\\
&\le&|\hat{\mu}(\hat{f}^{-j}(\hat{A})\cap \hat{B})-\mu(f^{-j}(\pi_N(\hat{A}))\cap \pi_N(\hat{B}))|\\
& &+|\mu(\pi_N(\hat{A}))\mu(\pi_N(\hat{B}))-\hat{\mu}(\hat{A})\hat{\mu}(\hat{B})|\\
&<&\epsilon+|\mu(\pi_N(\hat{A}))(\mu(\pi_N(\hat{B}))-\hat{\mu}(\hat{B}))|+|(\mu(\pi_N(\hat{A}))-\hat{\mu}(\hat{A}))\hat{\mu}(\hat{B})|\\
&<&\epsilon+\epsilon(\hat{\mu}(\hat{A})+\epsilon)+\epsilon\hat{\mu}(\hat{B})\\
&=&\epsilon(1+\hat{\mu}(\hat{A})+\hat{\mu}(\hat{B}))+\epsilon^2\\
&<&4\epsilon.
\end{eqnarray*}

If $f$ is ergodic, then 
$$\lim\limits_{n\rightarrow\infty}\frac{1}{n}\sum\limits_{j=0}^{n-1}\mu(f^{-j}(\pi_N(\hat{A}))\cap \pi_N(\hat{B}))=\mu(\pi_N(\hat{A}))\mu(\pi_N(\hat{B}))$$
or equivalently
$$\lim\limits_{n\rightarrow\infty}\frac{1}{n}\sum\limits_{j=0}^{n-1}\mu(f^{-j}(\pi_N(\hat{A}))\cap \pi_N(\hat{B}))-\mu(\pi_N(\hat{A}))\mu(\pi_N(\hat{B}))=0.$$
Therefore
$$|\lim\limits_{n\rightarrow\infty}\frac{1}{n}\sum\limits_{j=0}^{n-1}\hat{\mu}(\hat{f}^{-j}(\hat{A})\cap \hat{B})-\hat{\mu}(\hat{A})\hat{\mu}(\hat{B})|<4\epsilon.$$
Since $\epsilon$ is arbitrary, we must have
$\lim\limits_{n\rightarrow\infty}\frac{1}{n}\sum\limits_{j=0}^{n-1}\hat{\mu}(\hat{f}^{-j}(\hat{A})\cap \hat{B})-\hat{\mu}(\hat{A})\hat{\mu}(\hat{B})=0$,
i.e. $\hat{f}$ is ergodic.

If $f$ is weak-mixing, then
$$\lim\limits_{n\rightarrow\infty}\frac{1}{n}\sum\limits_{j=0}^{n-1}|\mu(f^{-j}(\pi_N(\hat{A}))\cap \pi_N(\hat{B}))-\mu(\pi_N(\hat{A}))\mu(\pi_N(\hat{B}))|=0.$$
Therefore
$$\lim\limits_{n\rightarrow\infty}\frac{1}{n}\sum\limits_{j=0}^{n-1}|\hat{\mu}(\hat{f}^{-j}(\hat{A})\cap \hat{B})-\hat{\mu}(\hat{A})\hat{\mu}(\hat{B})|<4\epsilon.$$
Since $\epsilon$ is arbitrary, we must have
$\lim\limits_{n\rightarrow\infty}\frac{1}{n}\sum\limits_{j=0}^{n-1}|\hat{\mu}(\hat{f}^{-j}(\hat{A})\cap \hat{B})-\hat{\mu}(\hat{A})\hat{\mu}(\hat{B})|=0$,
i.e. $\hat{f}$ is weak-mixing.

If $f$ is mixing, then 
$$\lim\limits_{n\rightarrow\infty}\mu(f^{-n}(\pi_N(\hat{A}))\cap \pi_N(\hat{B}))=\mu(\pi_N(\hat{A}))\mu(\pi_N(\hat{B}))$$
or equivalently
$$\lim\limits_{n\rightarrow\infty}\mu(f^{-n}(\pi_N(\hat{A}))\cap \pi_N(\hat{B}))-\mu(\pi_N(\hat{A}))\mu(\pi_N(\hat{B}))=0.$$
Therefore
$$|\lim\limits_{n\rightarrow\infty}\hat{\mu}(\hat{f}^{-n}(\hat{A})\cap \hat{B})-\hat{\mu}(\hat{A})\hat{\mu}(\hat{B})|<4\epsilon.$$
Since $\epsilon$ is arbitrary, we must have $\lim\limits_{n\rightarrow\infty}\hat{\mu}(\hat{f}^{-n}(\hat{A})\cap \hat{B})-\hat{\mu}(\hat{A})\hat{\mu}(\hat{B})=0$,
i.e. $\hat{f}$ is mixing.
\end{proof}

We now study the entropy properties of $\hat{f}$.

\begin{thm}\label{T:Entropy}
The topological entropy of $\hat{f}$ on $\hat{M}$ is equal to the topological entropy of $f$ on $M$.
\end{thm}
\begin{proof}
Since the canonical projection semiconjugates $\hat{f}$ to $f$, we have $h(\hat{f})\ge h(f)$ (see \cite{KH:Intro}, p.111).
Hence we only need to show that $h(\hat{f})\le h(f)$.

Fix $\epsilon>0$, and choose $l\in \mathbf{N}$ such that $2^{-l}\le \epsilon$. Let $Y$ be a $(n+l,\epsilon)$-spanning set
for $M$ with respect to $f$. For any $\hat{x}\in \hat{M}$, there exists $y\in Y$ such that 
$d_{n+l}^f(x_{-l},y)\le \epsilon$. Let $\hat{y}$ be a history such that $y_{-l}=y$. We have
\begin{eqnarray*}
\hat{d}_n^{\hat{f}}(\hat{x},\hat{y})&=&\max\limits_{0\le j\le n-1}\ \hat{d}(\hat{f}^j(\hat{x}),\hat{f}^j(\hat{y}))\\
&=&\max\limits_{0\le j\le n-1}\ \sum\limits_{i=0}^{\infty}2^{-i}\ d(f^j(x_{-i}),f^j(y_{-i}))\\
&=&\max\limits_{0\le j\le n-1}\ (\sum\limits_{i=0}^{l-1}2^{-i}\ d(f^j(x_{-i}),f^j(y_{-i})))+\sum\limits_{i=l}^{\infty}2^{-i}\ d(f^j(x_{-i}),f^j(y_{-i})))\\
&\le&\sum\limits_{i=0}^{l-1}2^{-i}\ \max\limits_{0\le j\le n-1}\ d(f^j(x_{-i}),f^j(y_{-i})))+\sum\limits_{i=l}^{\infty}2^{-i+l}2^{-l}\\
&\le&\sum\limits_{i=0}^{l-1}2^{-i}\epsilon+\sum\limits_{i=l}^{\infty}2^{-i+l}\epsilon\\
&<&4\epsilon.
\end{eqnarray*}
Hence $r_n^{\hat{f}}(4\epsilon,\hat{M})\le r_{n+l}^f(\epsilon,M)$. Therefore $h(\hat{f})\le h(f)$.
\end{proof}

\begin{thm}\label{T:Entropy1}
The measure-theoretic entropy $h_{\hat{\mu}}(\hat{f})$ of $\hat{f}$ with respect to the history measure $\hat{\mu}$ 
is equal to the measure-theoretic entropy $h_{\mu}(f)$ of $f$ with respect to $\mu$.
\end{thm}
\begin{proof}
For any $\epsilon>0$, $n\in \mathbf{N}$ and $\hat{x}\in \pi_0^{-1}(x)$, if $\hat{y}\in B_n^{\hat{f}}(\hat{x},\epsilon)$
then $y_0\in B_n^f(x,\epsilon)$ since $\hat{d}(\hat{x},\hat{y})=\sum_{i=0}^{\infty}2^{-i}\ d(x_{-i},y_{-i})$. Hence,
$\hat{\mu}(B_n^{\hat{f}}(\hat{x},\epsilon))\le \mu(B_n^f(x,\epsilon))$. This implies that
$h_{\hat{\mu}}(\hat{f},\hat{x})\ge h_{\mu}(f,x)$ for $\mu$-almost every $x\in M$, by Theorem \ref{T:BK} (a).

On the other hand, for any $\delta>0$ we can choose $\epsilon>0$ such that 
$$\limsup\limits_{n\rightarrow\infty}-\frac{1}{n}\ \log\ \hat{\mu}(B_n^{\hat{f}}(\hat{x},4\epsilon))\ge h_{\hat{\mu}}(\hat{f},\hat{x})-\delta.$$
Choose $j\ge 0$ such that $2^{-j}\le \epsilon$. We also have
$$\limsup\limits_{n\rightarrow\infty}-\frac{1}{n}\ \log\ \mu(B_{n+j}^f(x_j,\epsilon))\le h_{\mu}(f,x_j)=h_{\mu}(f,x),$$
where the inequality follows from Theorem \ref{T:BK} (a) and the equality follows from Theorem \ref{T:BK} (b).
By the estimate we get in Theorem \ref{T:Entropy}, we have that if $y\in B_{n+j}^f(x_j,\epsilon)$ then 
$\hat{y}\in B_n^{\hat{f}}(\hat{x},4\epsilon)$ for any $\hat{y}\in \pi_j^{-1}(y)$. Hence,
$\hat{\mu}(B_n^{\hat{f}}(\hat{x},4\epsilon))\ge \mu(B_{n+j}^f(x_j,\epsilon))$. Letting $n$ go to infinity and
using the two inequalities above, we get $h_{\hat{\mu}}(\hat{f},\hat{x})\le h_{\mu}(f,x)+\delta$. Since
$\delta$ is arbitrary, we have that $h_{\hat{\mu}}(\hat{f},\hat{x})\le h_{\mu}(f,x)$, which gives us
$h_{\hat{\mu}}(\hat{f},\hat{x})=h_{\mu}(f,x)$ for $\mu$-almost every $x\in M$ and $\hat{\mu}$-almost every $\hat{x}\in \pi_0^{-1}(x)$. 
Therefore by Theorem \ref{T:BK} (c)
$$h_{\hat{\mu}}(\hat{f})=\int\limits_{\hat{M}}\ h_{\hat{\mu}}(\hat{f},\hat{x})\ d\hat{\mu}
=\int\limits_{\pi_0^{-1}(M)}\ h_{\hat{\mu}}(\hat{f},\hat{x})\ d\hat{\mu}
=\int\limits_{M}\ h_{\mu}(f,x)\ d\mu =h_{\mu}(f).$$
\end{proof}

As a corollary to the above two theorems,  we have the following.

\begin{cor}\label{C:Entropy}
The history measure $\hat{\mu}$ is a measure of maximal entropy of $\hat{f}$  if and only if $\mu$ is a measure 
of maximal entropy of $f$. And $\hat{\mu}$ is the unique measure of maximal entropy of $\hat{f}$ if and only if 
$\mu$ is the unique measure of maximal entropy of $f$.
\end{cor}
\begin{proof}
By Theorem \ref{T:Entropy}, the topological entropy of $\hat{f}$ on $\hat{M}$ is equal to the topological entropy 
of $f$ on $M$. And by Theorem \ref{T:Entropy1}, the measure-theoretic entropy of $\hat{f}$ with respect to the 
history measure $\hat{\mu}$ is equal to the measure-theoretic entropy of $f$ with respect to $\mu$. Hence, 
$\hat{\mu}$ is a measure of maximal entropy of $\hat{f}$  if and only if $\mu$ is a measure of maximal entropy of $f$.

If $\mu$ and $\nu$ are two different measures of maximal entropy of $f$, then there exists a set $A\subset M$
such that $\mu(A)\neq \nu(A)$. Then, $\hat{\mu}(\pi_0^{-1}(A))\neq \hat{\nu}(\pi_0^{-1}(A))$. So $\hat{\mu}$
and $\hat{\nu}$ are two different measures of maximal entropy of $\hat{f}$. Similarly, if $\hat{\mu}$ and
$\hat{\nu}$ are two different measures of maximal entropy of $\hat{f}$, then there exists a set $\hat{A}\subset \hat{M}$
such that $\hat{\mu}(\hat{A})\neq \hat{\nu}(\hat{A})$. By the definition of the history measure, there exists
$j\ge 0$ such that $\mu(\pi_j(\hat{A}))\neq \nu(\pi_j(\hat{A}))$. So $\mu$ and $\nu$ are two different measures 
of maximal entropy of $f$.
\end{proof}

Periodic points for $\hat{f}$ and periodic points for $f$ are also closely related. Recall that $E_n(f)$ is the set of periodic
points for $f$ whose periods are factors of $n$ and $N_n(f)$ is the cardinality of $E_n(f)$.

\begin{lem}\label{L:Periodic}
(a) For any $j\ge 0$, $\pi_j(E_n(\hat{f}))=E_n(f)$. In particular $N_n(\hat{f})=N_n(f)$;
(b) Periodic points for $\hat{f}$ are dense in $\hat{M}$ if and only if periodic points for $f$ are dense in $M$.
\end{lem}
\begin{proof}
If $\hat{x}\in E_n(\hat{f})$, then $\hat{f}^n(\hat{x})=\hat{x}$, i.e. $f^n(x_{-i})=x_{-i}$ for any $i\ge 0$. On the
other hand, if $x\in E_n(f)$, then $f^n(x)=x$ and $f^n(f^j(x))=f^j(x)$ for any $j\ge 0$. So $\hat{x}$ with
$x_{-j}=f^{n-j}(x)$ for $0\le j < n$ and $x_{-i-n}=x_{-i}$ for any $i\ge 0$ is in $E_n(\hat{f})$. Hence $\pi_j$ induces a bijection 
$\tilde{\pi}_j$ between $E_n(\hat{f})$ and $E_n(f)$ for any $j\ge 0$. This proves part (a).

Assume that periodic points for $\hat{f}$ are dense in $\hat{M}$. Let $U\subset M$ be any
open set. By the assumption, there is a periodic point $\hat{x}$ for $\hat{f}$ in $\pi_0^{-1}(U)$, which means
that $x_0$ is a periodic point for $f$ in $U$. Since $U$ is arbitrary, we have that periodic points for $f$ are
dense in $M$. Conversely, assume that periodic points for $f$ are dense in $M$. Let
$\hat{U}\subset \hat{M}$ be any open set. Without loss of generality, we can write $\hat{U}=\pi_j^{-1}(V)$ with $V$ open in $M$ 
for some $j\ge 0$. By the assumption, there is a periodic point $x$ for $f$ in $V$, which means that $\tilde{\pi}_j^{-1}(x)$ is a 
periodic point for $\hat{f}$ in $\pi_j^{-1}(V)$. Since $\hat{U}$ is arbitrary, we have that periodic points for $\hat{f}$ are dense in $\hat{M}$.
\end{proof}

\begin{thm}\label{T:Periodic}
The history measure $\hat{\mu}$ describes the distribution of periodic points for $\hat{f}$ if and only if
$\mu$ describes the distribution of periodic points for $f$.
\end{thm}
\begin{proof}
By Lemma \ref{L:Periodic}(a), we have 
$$(\pi_j)_{*}(\frac{1}{N_n(\hat{f})}\sum\limits_{\hat{x}\in E_n(\hat{f})}\ \delta_{\hat{x}})=
\frac{1}{N_n(f)}\sum\limits_{x\in E_n(f)}\ \delta_x,\ \ \forall j\ge 0.$$
By the definition of the history measure $\hat{\mu}$, it follows that 
$$\frac{1}{N_n(\hat{f})}\sum\limits_{\hat{x}\in E_n(\hat{f})}\ \delta_{\hat{x}}\ \longrightarrow \hat{\mu}$$
if and only if 
$$\frac{1}{N_n(f)}\sum\limits_{x\in E_n(f)}\ \delta_x \longrightarrow \mu.$$
\end{proof}

Finally, we make an easy observation concerning Lyapunov exponents. For detailed information on
this very important concept, see \cite{KH:Intro}.

Assume that $f$ is differentiable on $M$. The pullback under $\pi_0$ of the tangent bundle $T_M$ of $M$ is 
a bundle on $\hat{M}$, the tangent bundle $T_{\hat{M}}$. The derivative $Df$ of $T_M$ lifts to a map $D\hat{f}$ of $T_{\hat{M}}$,
the derivative of $\hat{f}$, in a natural way. Note that $\hat{f}$ is invertible on the history space $\hat{M}$, therefore the map
$D\hat{f}^{-1}$ is also well defined. A point in $T_{\hat{M}}$ is of the form $(\hat{x}, \hat{v})$, where $\hat{x}\in \hat{M}$ and 
$\hat{v}=\{v_{-i}\}_{i=0}^{\infty}$ is such that $v_0\in T_M(x_0)$ and $v_{-i}=D_{x_0}\hat{f}^{-i}(v_0)$.
Since $\hat{v}$ is uniquely determined by $\hat{x}$ and $v_0$, we can also denote a point in $T_{\hat{M}}$ by
$(\hat{x}, v)$ with $\hat{x}\in \hat{M}$ and $v\in T_M(x_0)$. 

The Lyapunov exponent for $f$ at a point $(x, v)\in T_M$ is defined as the limit $\lim_{n\rightarrow \infty} \frac{1}{n} \| D_xf^n(v) \|$,
if it exists, where $\|\cdot\|$ denotes a norm on the tangent bundle. By the above discussion, we have the following.

\begin{lem}\label{L:Lyapunov}
The Lyapunov exponent for $\hat{f}$ at $(\hat{x}, v)\in T_{\hat{M}}$ is equal to the Lyapunov exponent for $f$
at $(x_0, v)\in T_M$.
\end{lem}
\begin{proof}
Note that, by definition, Lyapunov exponents are invariant under the map $(f, Df)$.
\end{proof}

\section{Nonalgebraic Attractors}\label{S:Attractor}

The goal of this section is to prove Theorem \ref{T:Main1}.

Recall that $f_\lambda$ is a holomorphic map on $\mathbf{P}^k$, $k\ge2$, of the following form 
$$f_\lambda[z:w_1:\cdots:w_{k-1}:t] = [(z-2w_1)^2:\cdots:(z-2w_{k-1})^2:z^2:t^2+\lambda z^2].$$

Let us start with investigating the dynamics of $f_0$ on the set $\Pi = \{t = 0\}$, the hyperplane at infinity. 

\subsection{The Map $f=f_0|_{\Pi}$}

Before we study the properties of the map $f$, we need the following definition and result from the study of critically finite maps (see \cite{Ro:Critical} for details).

\begin{dfn}
Let $g:\mathbf{P}^n\rightarrow \mathbf{P}^n$ be a holomorphic map. Let $C_1$ be the critical set of $g$ given by
$$C_1 = \{ p\in \mathbf{P}^k | rank(dg(p)) < k-1\},$$
where $dg(p)$ denotes the differential of $g$ at $p$. Let $D_1$ be the post-critical set given by
$$D_1 = \bigcup\limits_{j=1}^{\infty} f^j(C_1),$$
and $E_1$ be the $\omega-$limit set given by
$$E_1 = \bigcap\limits_{j=1}^{\infty} f^j(\overline{D_1}).$$
We say that $g$ is critically finite if $D_1$ is an analytic (hence algebraic) set in $\mathbf{P}^n$.
In this case, we inductively define for $j\in\{2,\cdots,n\}$
$$C_j = C_1\bigcap E_{j-1},$$
$$D_j = \bigcup\limits_{j=1}^{\infty}f^j(C_j),$$
$$E_j = \bigcap\limits_{j=1}^{\infty}f^j(\overline{D_j}).$$
We say that $g$ is $n-$critically finite if $D_j$ is algebraic for $1\le j\le n$ and none of the irreducible components of $E_n$ are contained in $C_1$.
\end{dfn}

\begin{thm}\label{T:Critical}
Let $g:\mathbf{P}^n\rightarrow \mathbf{P}^n$ be a holomorphic map. If $g$ is $n-$critically finite, then the Julia set $J_1$ of $g$ is the whole of $\mathbf{P}^n$.
\end{thm}

For the map $f[z_0:z_1:\cdots:z_{k-1}] = [(z_0-2z_1)^2:\cdots:(z_0-2z_{k-1})^2:z_0^2]$, we have the following.

\begin{prop}\label{P:Main}
The holomorphic map $f$ is $(k-1)-$critically finite and $\mathbf{P}^{k-1}$ is the only nonempty closed backward invariant subset of $\mathbf{P}^{k-1}$.
\end{prop}

\begin{rmk}\label{R:Julia}
Note that all the Julia sets $J_m,1\le m\le k-1$, are nonempty closed backward invariant subsets of $\mathbf{P}^{k-1}$. Therefore they are 
all equal to $\mathbf{P}^{k-1}$ for the given map $f$.
\end{rmk}

The following theorem is a generalization and improvement of Proposition 7.5 in \cite{FS:DynamicsI} and we give
two different proofs, one analytic and the other algebraic.

\begin{thm}\label{T:Injective}
Let $g:\mathbf{P}^n\rightarrow \mathbf{P}^n$ be a holomorphic map of algebraic degree $d>1$. Let $X\subset \mathbf{P}^n$ be an 
irreducible algebraic set of pure dimension $s$ and assume that $g(X)=X$. Then the topological degree $l$ of $g$ restricted to $X$
is equal to $d^s$.
\end{thm}
\begin{proof}[Proof 1 (Analytic)]
Let $\omega$ be the standard K$\ddot{a}$hler form on $\mathbf{P}^n$. Since the pullback operator $g^*$ commutes with the
differential operator $d$ on forms, $g^*$ induces an operator on the cohomology groups of $\mathbf{P}^n$, also denoted by
$g^*$. Hence we have
$$[g^*(\omega)]=d[\omega],\ \ \ \ \ \ [g^*(\omega^s)]=d^s[\omega^s].$$
Therefore $g^*(\frac{\omega^s}{s!})-d^s\frac{\omega^s}{s!}=d\varphi$, where $\varphi=\varphi_1+\varphi_2$ with $\varphi_1$ a form of
bidegree $(s-1,s)$ and $\varphi_2$ a form of bidegree $(s,s-1)$. (In fact, we can write $d\varphi=dd^c\varphi^\prime$,
where $\varphi^\prime$ is a form of bidegree $(s-1,s-1)$. This is a consequence of Hodge theory on $\mathbf{P}^b$.)

Let $[X]$ denote the current of integration on $X$. We have
$$\int[X]\wedge g^*(\frac{\omega^s}{s!}) = d^s\int[X]\wedge \frac{\omega^s}{s!} + \int[X]\wedge d\varphi.$$
The left-hand side is equal to $l\cdot volume(X)$, since $g(X)=X$ and $l$ is the topological degree of $g$ restricted
to $X$. The first term on the right-hand side is equal to $d^s\cdot volume(X)$ by Wirtinger's theorem and the second
term vanishes by Stoke's theorem. Hence we have $l\cdot volume(X)=d^s\cdot volume(X)$, i.e. $l=d^s$.
\end{proof}
\begin{proof}[Proof 2 (Algebraic)]
Assume $X$ is an algebraic set of degree $r$, i.e. for a generic choice of $(n-s)-$dimensional plane $H$ in $\mathbf{P}^n$,
$X\cdot H=r$. Without loss of generality, we can choose $H$ to be the plane $\{z_0=0\}\cap\cdots \cap\{z_{s-1}=0\}$. 
Write $g=[g_0:\cdots:g_n]$ and $z=[z_0:\cdots:z_n]$, then $g^{-1}(H)=\{g_0(z)=0\}\cap\cdots \cap\{g_{s-1}(z)=0\}$. 
Since $g(X)=X$ and $l$ is the topological degree of $g$ restricted to $X$, $X\cdot g^{-1}(H)=l\cdot r$. On the other hand, 
each $\{g_j(z)=0\}$ , $0\le j<s$, is linearly equivalent to $d\cdot E$, where $E$ is the hyperplane divisor at infinity (see \cite{GH:Principle}). 
Hence $g^{-1}(H)$ is linearly equivalent to $d^s\cdot E^s$. Since $X\cdot E^s=r$, we also have $X\cdot f^{-1}(H)=d^s\cdot r$.
Therefore, $l=d^s$.
\end{proof}

\begin{rmk}
If $X$ is biholomorphic to $\mathbf{P}^s$, then $g$ restricted to $X$ actually has algebraic degree $d$. This follows from the fact that 
every holomorphic map on $\mathbf{P}^s$ is given by $s+1$ homogeneous polynomials (\cite[Theorem 2.1]{FS:DynamicsI}).
\end{rmk}

\begin{lem}\label{L:Subcritical}
Let $g:\mathbf{P}^n\rightarrow \mathbf{P}^n$ be a $n-$critically finite map of degree $d>1$. Let $X\subset E_m,\ 1\le m< n$,
be an irreducible component of $E_m$, which is fixed by $g^l$ for some $l\ge 1$. If $X$ is biholomorphic to $\mathbf{P}^{n-m}$,
then $g^l$ restricted to $X$ is a $(n-m)-$critically finite map of algebraic degree $d^l$.
\end{lem}
\begin{proof}
Denote $\tilde{g}=g^l|_X$. A point $p\in X$ is critical for $\tilde{g}$ if and only if $\tilde{g}$ is
non-injective in $U\cap X$ for every neighborhood $U$ of $p$ in $\mathbf{P}^k$, which can only happen if $g^l$ is 
non-injective in $U$. The latter is only possible if $p$ is a critical point for $g^l$. Hence, the critical set 
$\tilde{C_1}$ for $\tilde{g}$ is contained in $\cup_{i=0}^{l-1}g^{-i}(C_1)\cap X$. Since every point $q\in g^{-i}(C_1)\cap X$, $0\le i\le l-1$,
is mapped to a point $g^i(q)\in C_1\cap E_m$, the post-critical set $\tilde{D_1}$ for $\tilde{g}$ is contained in
$D_{m+1}$, and the $\omega-$limit set $\tilde{E_1}$ for $\tilde{g}$ is contained in $E_{m+1}$. Then by the same argument
we can show that $\tilde{E_j}\subset E_{m+j}$, $2\le j\le n-m$. If $Y$ is a component of $\tilde{E_j}$ that is contained
in $\tilde{C_j}$, $1\le j\le n-m$, then its image under some iterations of $g$ will be a component of $E_{m+j}$ contained in
$C_{m+j}$, which is impossible since $g$ is $n-$critically finite. It is easily seen that $\tilde{D_j}$ is algebraic for $1\le j\le n-m$, hence
$\tilde{g}$ is $(n-m)-$critically finite. The fact that $\tilde{g}$ has algebraic degree $d^l$ follows from the remark above.
\end{proof}

The following proposition is due to M. Jonsson.

\begin{prop}\cite[Proposition 3.6]{J:Critical}\label{P:Jonsson}
If $g:\mathbf{P}^n\rightarrow \mathbf{P}^n$ is $n-$critically finite, then any nonempty closed backward invariant subset $E\subset \mathbf{P}^n$
must contain all the repelling periodic points for $g$.
\end{prop}

\begin{proof}[Proof of Proposition \ref{P:Main}]
Obviously, the critical set $C_1$ of $f$ is $\{z_0=0\}\cup(\cup_{j=1}^{k-1}\{z_0=2z_j\})$ and the $\omega-$limit set
$E_1$ of $f$ is $\cup\{z_i=z_j\}$, $0\le i < j\le k-1$. Consequently, the set $E_m$, $2\le m\le k-1$, contains intersections
of $m$ different hyperplanes in $E_1$. In particular, $E_{k-1}$ consists of the single point $[1:\cdots:1]$, which is
a repelling fixed point. Hence $f$ is $(k-1)-$critically finite.

Let $V_1$ denote the set of preimages of $E_1$. As in \cite{FS:Critical}, we show that $\mathbf{P}^{k-1}\backslash V_1$ is hyperbolic. $V_1$ contains the 
set $H$ of $2k-1$ hyperplanes $\{z_i=0\}$, $0\le i\le k-1$, and $\{z_0=2z_j\}$, $1\le j\le k-1$. We can use coordinates so that $\{z_0=0\}$ is the hyperplane 
at infinity. That is, we set the homogeneous coordinate $z_0=1$ and consider the remaining $\mathbf{C}^{k-1}$ with coordinates $(z_1,\cdots,z_{k-1})$.
Then $V_1$ contains the affine complex hyperplanes $\{z_j=0\}$ and $\{z_j=1/2\}$, $1\le j\le k-1$. Hence $\mathbf{P}^{k-1}\backslash V_1$ is contained in
a product of $k-1$ copies of $\bar{\mathbf{C}}\backslash \{0,1/2,\infty\}$. The universal covering of this set is the unit polydisc in $\mathbf{C}^{k-1}$
and is therefore Kobayashi hyperbolic. Since holomorphic maps and inclusion maps are distance decreasing $\mathbf{P}^{k-1}\backslash V_1$ is
also hyperbolic. Let $U$ be an open ball in $\mathbf{P}^{k-1}\backslash V_1$ such that the restriction of the family of iterates $\{f^i\}$ to $U$ is a family 
of holomorphic maps of $U$ into $\mathbf{P}^{k-1}\backslash V_1$. Then this family is normal, which contradicts Theorem \ref{T:Critical}. 
Hence $V_1$ is dense in $\mathbf{P}^{k-1}$.

Similarly, the set $V_m$ of preimages of $E_m$ contains the set $H\cap E_{m-1}$, $2\le m\le k$. Since every irreducible
component $X$ of $E_{m-1}$ is the intersection of $m-1$ different hyperplanes in $E_1$, $X\backslash V_m$ is contained in a product
of $k-m$ copies of $\bar{\mathbf{C}}\backslash \{0,1/2,\infty\}$, hence hyperbolic. By Lemma \ref{L:Subcritical}, if $f^l$ maps $X$ into itself, then
$f^l$ restricted to $X$ is $(k-m)-$critically finite. Then by Theorem \ref{T:Critical}, we conclude that $V_m\cap X$ is dense in 
$X$ as above. Hence $V_m$ is dense in $E_{m-1}$. 

Let $B$ be a nonempty closed backward invariant subset of $\mathbf{P}^{k-1}$. By Proposition \ref{P:Jonsson}, the set $B$ contains the point $[1:\cdots:1]$,
i.e. the set $E_k$, therefore a dense set in $\mathbf{P}^{k-1}$. But $B$ is closed, so $B=\mathbf{P}^{k-1}$.
\end{proof}

From the above discussion, we know that the support of the Green measure associated with $f$ is equal to $\Pi$. We need the following result 
from \cite{FS:Potential}.

\begin{lem}\label{L:FS}
Let $g$ be a holomorphic map on $\mathbf{P}^n$ and $\mu$ be the associated Green measure on
$\mathbf{P}^n$. Let $U$ be an open set intersecting the support of $\mu$. Define
$E_U:=\mathbf{P}^n\backslash\cup_{N\ge 0}(\cap_{j\ge N}g^j(U))$. Then $E_U$ is locally pluripolar in $\mathbf{P}^n$.
\end{lem}

\begin{cor}\label{C:Mixing}
Let $g$ be a holomorphic map on $\mathbf{P}^n$ and $\mu$ be the associated Green measure on
$\mathbf{P}^n$. If the support of $\mu$ is equal to $\mathbf{P}^n$, then $g$ is topologically mixing on
$\mathbf{P}^n$. 
\end{cor}
\begin{proof}
Assume that $g$ is not topologically mixing on $\mathbf{P}^n$. Then by definition, there exist two
open sets $\Omega_1$ and $\Omega_2$ in $\mathbf{P}^n$ and arbitrarily large $j$ such that 
$g^j(\Omega_1)\cap \Omega_2 = \emptyset$. But then we have
$E_{\Omega_1} = \mathbf{P}^n\backslash\cup_{N\ge 0}(\cap_{j\ge N}g^j(\Omega_1))\supset \Omega_2$, which
contradicts Lemma \ref{L:FS} and proves the corollary.
\end{proof}

As an immediate corollary, we have

\begin{cor}\label{C:Mixing1}
The map $f$ is topologically mixing on $\Pi$.
\end{cor}

\subsection{Attractors}

We now study the properties of $f_\lambda$.

The following lemma gives a "trapping region" near $\Pi$ for $f_\lambda$.

\begin{lem}\label{L:Trapping}
For $\rho > 0$, let $U_\rho\subset \mathbf{P}^k$ be the neighborhood of $\Pi$ defined by 
$$U_\rho = \{[z:w_1:\cdots:w_{k-1}:t]: |t|<\rho |(z,w_1,\cdots,w_{k-1})|\},$$
where $|(z,w_1,\cdots,w_{k-1})| = \max\{|z|,|w_1|,\cdots,|w_{k-1}|\}$. If $0 < 2 |\lambda| < \rho < \sqrt{|\lambda|}$
then $f_\lambda(U_\rho)\subset\subset U_\rho$.
\end{lem}
\begin{proof}
For any $[z:w_1:\cdots:w_{k-1}:t] \in U_\rho$,
write $f_\lambda [z:w_1:\cdots:w_{k-1}:t] = [z^\prime:w_1^\prime:\cdots:w_{k-1}^\prime:t^\prime]$.

If $|z|\ge |w_j|$ for all $j\in\{1,\cdots,k-1\}$, then $|t^\prime| = |t^2 + \lambda z^2|\le 
|t|^2 + |\lambda||z|^2 < (\rho^2 + |\lambda|)|z|^2 < 2|\lambda||z|^2 = 2|\lambda||w_{k-1}^\prime|
\le 2|\lambda||(z^\prime,w_1^\prime,\cdots,w_{k-1}^\prime)| < 
\rho|(z^\prime,w_1^\prime,\cdots,w_{k-1}^\prime)|$.

If $|z| < |w_j|$ for some $j\in\{1,\cdots,k-1\}$ we assume without loss of generality that 
$|w_1| = \max\{|w_j|: |w_j| > |z|, j\in\{1,\cdots,k-1\}\}$. Then $|t^\prime| = |t^2 + \lambda z^2|
\le |t|^2 + |\lambda||z|^2 < (\rho^2 + |\lambda|)|w_1|^2 < 2|\lambda| |z^\prime|
\le 2|\lambda||(z^\prime,w_1^\prime,\cdots,w_{k-1}^\prime)| <
\rho|(z^\prime,w_1^\prime,\cdots,w_{k-1}^\prime)|$.
\end{proof}

The following semiconjugation property is a key observation for our study (\cite{JW:Nonalgebraic}). We include a proof for the
convenience of the reader. 

\begin{prop}\label{P:Semi}
Let $U_\rho$ be the trapping region defined in Lemma \ref{L:Trapping}, with $0 < 2 |\lambda| < \rho < \sqrt{|\lambda|}$. If
$$K_\lambda := \bigcap_{n\ge 0} f_\lambda^n(U_\rho),$$
then $K_\lambda\subset U_\rho$ and $f_\lambda(K_\lambda) = K_\lambda$. Moreover, there exists a continuous and onto 
map $\phi_\lambda$ between $\hat{\Pi}$ and $K_\lambda$ such that the following diagram commutes, where 
$\pi[z:w_1:\cdots:w_{k-1}:t] = [z:w_1:\cdots:w_{k-1}]$ and $\hat{f_0}$ is the natural extension of $f_0$ on $\hat{\Pi}$, 
the history space associated with $f_0$ on $\Pi$ .
\[
  \begin{CD}
    \hat{\Pi} @>\hat{f_0}>> \hat{\Pi}\\
    @V\phi_\lambda VV    @VV\phi_\lambda V\\
    K_\lambda @>f_\lambda >> K_\lambda\\
    @V\pi VV   @VV\pi V\\
    \Pi @>f_0>> \Pi
  \end{CD}
\]  
\end{prop}
\begin{proof}
The first two assertions follow easily from Lemma \ref{L:Trapping}. We prove the last statement by explicitly constructing
the map $\phi_\lambda : \hat{\Pi} \rightarrow K_\lambda$.

For $a\in \Pi$, let $L_a = \overline{\pi^{-1}(a)}$ be the line in $\mathbf{P}^k$ passing through $a$ and $[0:0:\cdots:0:1]$. Note that $f_\lambda$ 
fixes $[0:0:\cdots:0:1]$ and maps lines to lines. So we have $f_\lambda(L_a) = L_{f(a)}$. If $V_a = U_\rho\cap L_a$, 
then $a\in V_a$ and Lemma \ref{L:Trapping} shows that $f_\lambda(V_a)\subset\subset V_{f(a)}$ for all $a\in \Pi$. 
Therefore the diameter of the set $f_\lambda^n(V_a)$ tends to zero, uniformly in $a$, as $n \to \infty$. 
Hence, for any $\hat{a}\in \hat{\Pi}$, the intersection $\bigcap_{i\ge 0} f_\lambda^i(V_{a_{-i}})$ is a single point, which we denote by $\phi_\lambda(\hat{a})$. 

We now show that $\phi_\lambda$ maps $\hat{\Pi}$ onto $K_\lambda$. For any $\hat{a}\in \hat{\Pi}$, we have $\phi_\lambda(\hat{a})\in f_\lambda^n(U_\rho)$ 
for all $n\ge 0$ so $\phi_\lambda(\hat{a})\in K_\lambda$. Conversely, if $x\in K_\lambda$, then $x$ has a prehistory $\hat{x} = \{x_{-i}\}$ in $K_\lambda$.
If $a_{-i} = \pi(x_{-i})$, then $\hat{a} = \{a_{-i}\}\in \hat{\Pi}$ and $x = \phi_\lambda(\hat{a})$.

The continuity of $\phi_\lambda$ and the semiconjugation properties of $\phi_\lambda$ and $\pi$ follow easily from the construction. 
\end{proof}

\begin{rmk}\label{R:Semi}
From the proof of the above proposition, we see that the map $\phi_\lambda:\hat{\Pi}\rightarrow K_\lambda$ can be lifted to a map
$\hat{\phi_\lambda}:\hat{\Pi}\rightarrow \hat{K_\lambda}$ such that $\phi_\lambda=\pi_\lambda\circ \hat{\phi_\lambda}$,
where $\hat{K_\lambda}$ is the history space associated with $f_\lambda|_{K_\lambda}$ and 
$\pi_\lambda: \hat{K_\lambda}\rightarrow K_\lambda$ is the natural projection. Hence we have the following commutative diagram.
\[
  \begin{CD}
    \hat{\Pi} @>\hat{\phi_\lambda}>> \hat{K_\lambda}\\
    @V\pi_{\Pi} VV    @VV\pi_\lambda V\\
    \Pi @<\pi << K_\lambda\\
   \end{CD}
\]  
It is also obvious from the construction of $\phi_\lambda$ that $\hat{\phi_\lambda}$ is continuous, one-to-one and onto, hence a conjugation
between $\hat{\Pi}$ and $\hat{K_\lambda}$.
\end{rmk}

As a corollary, we show that the set $K_\lambda$ given in the previous proposition is indeed an attractor for $f_\lambda$. 

\begin{cor}\label{C:Attractor}
The $f_\lambda-$invariant set $K_\lambda$ given in Proposition \ref{P:Semi} is an attractor for $f_\lambda$.
\end{cor}
\begin{proof}
We are going to show that $f_\lambda$ is topologically mixing on $K_\lambda$. The assertion then follows from Remark \ref{R:Attractor}.

If $U$ and $V$ are two open sets in $K_\lambda$, then $\phi_\lambda^{-1}(U)$ and $\phi_\lambda^{-1}(V)$ are two open sets
in $\hat{\Pi}$. From Theorem \ref{T:Transitive} and Corollary \ref{C:Mixing1} we know that $\hat{f_0}$ is topologically mixing on $\hat{\Pi}$. 
So there exists $N>0$, such that for any $n\ge N$, $\hat{f_0}^n(\phi_\lambda^{-1}(U))\cap \phi_\lambda^{-1}(V)\neq \emptyset$. Hence 
$f_{\lambda}^n(U)\cap V\neq \emptyset$, i.e. $f_\lambda$ is topologically mixing on $K_\lambda$.
\end{proof}

We will next show that the attractor $K_\lambda$ is nonalgebraic.

To see this, we look at $M := K_\lambda\cap W$, where $W$ is the complex surface ($\mathbf{P}^2$) defined by
$$W := \{[z:w_1:\cdots:w_{k-1}:t]: w_1 = w_2 = \cdots = w_{k-1}\}.$$
Observe that
$$\begin{array}{c}
f_\lambda(W) = \{[z:w_1:\cdots:w_{k-1}:t]: z = w_1 = \cdots = w_{k-2}\}, \\
f_\lambda^2(W) = \{[z:w_1:\cdots:w_{k-1}:t]: z = w_1 = \cdots = w_{k-3} = w_{k-1}\}, \\
\vdots \\
f_\lambda^{k-1}(W) = \{[z:w_1:\cdots:w_{k-1}:t]: z = w_1 = w_3 = \cdots = w_{k-1}\}, \\
f_\lambda^k(W) = \{[z:w_1:\cdots:w_{k-1}:t]: w_1 = w_2 = \cdots = w_{k-1}\}. \\
\end{array}$$
So $W$ is invariant under $f_\lambda^k$ and we can define $g_\lambda := f_\lambda^k|_W$. Note that if 
$k=2$ then $W$ is the whole of $\mathbf{P}^2$. 

Let $L\subset W$ be the line fixed by $f_\lambda$ defined by
$$L := \{[z:w_1:\cdots:w_{k-1}:t]: z = w_1 = \cdots = w_{k-1}\}.$$
Let $p_\lambda$ be the fixed point $[1:1:\cdots:1:t_\lambda]\in L$, 
where $t_\lambda = \frac{1-\sqrt{1-4\lambda}}{2}$.
Let $q_\lambda$ be the point $[1:1:\cdots:1:-t_\lambda]\in L$.
Note that $f_\lambda(q_\lambda) = p_\lambda$ and $g_\lambda(q_\lambda) = p_\lambda$.
Denote $U = U_\rho\cap W$, with $0 < 2 |\lambda| < \rho < \sqrt{|\lambda|}$.
We need two computational lemmas.

\begin{lem}\label{L:Line}
There is no line contained in $U$ that is fixed by $g_\lambda$, if $|\lambda|\neq 0$ is sufficiently small.
\end{lem}
\begin{proof}
By the observation before the lemma, we have
$$\begin{array}{c}
f_\lambda([z:w:\cdots:w:t]) = [(z-2w)^2:\cdots:(z-2w)^2:z^2:t^2+\lambda z^2], \\
f_\lambda^2([z:w:\cdots:w:t]) = [(z-2w)^4:\cdots:(z-2w)^4:((z-2w)^2-2z^2)^2: \\
(z-2w)^4:(t^2+\lambda z^2)^2+(z-2w)^4], \\
\vdots \\
f_\lambda^k([z:w:\cdots:w:t]) = [P_1(z,w):(z-2w)^{2^k}:\cdots:(z-2w)^{2^k}:Q_1(z,w,t)], \\
\end{array}$$
where $$P_1(z,w) = ((z-2w)^{2^{k-1}}-\cdots-2((z-2w)^4-2((z-2w)^2-2z^2)^2)^2\cdots)^2,$$
$$Q_1(z,w,t) = (\cdots((t^2+\lambda z^2)^2 + \lambda(z-2w)^4)^2 +\cdots + \lambda(z-2w)^{2^{k-1}})^2 + \lambda(z-2w)^{2^k}.$$

Therefore, we can write $g_\lambda$ as
$$g_\lambda([z:1:\cdots:1:t]) = [P(z):1:\cdots:1:Q(z,t)],$$
where $$P(z) = (1-2(1-\cdots-2(1-2(1-\frac{2z^2}{(z-2)^2})^2)^2\cdots)^2)^2,$$
$$Q(z,t) = (\cdots(\frac{(t^2+\lambda z^2)^2}{(z-2)^4} + \lambda)^2 +\cdots + \lambda)^2 + \lambda.$$

Obviously, for all $\alpha,\beta\in \mathbf{C}$, the lines $\{t = \alpha\}$ are not fixed, and the lines
$\{z = \beta\}$ are not in $U$. So we only need to consider lines of the form $\{t = \alpha z +\beta\}$.

Let $z = \infty$, we get
\begin{equation}
|\alpha| < \rho\ and\ (\cdots((\alpha^2 + \lambda)^2 + \lambda)^2 +\cdots + \lambda)^2 + \lambda = \alpha + \beta.
\end{equation}
Let $z = 1$, we get
\begin{equation}
|\alpha + \beta| < \rho\ and\ (\cdots(((\alpha + \beta)^2 + \lambda)^2 + \lambda)^2 +\cdots + \lambda)^2 + \lambda = \alpha + \beta.
\end{equation}
Let $z = 0$, we get
\begin{equation}
|\beta| < \rho\ and\ (\cdots((\frac{\beta}{2})^4 + \lambda)^2 +\cdots + \lambda)^2 + \lambda = \alpha + \beta.
\end{equation}

From $(1)$ and $(2)$, we have $\beta = 0$ or $\beta = -2\alpha$. If $\beta = 0$, 
we have $\alpha = t_\lambda$ from $(1)$, which does not satisfy $(3)$. If $\beta = -2\alpha$,  
we have $\alpha = -t_\lambda$ from $(1)$, which does not satisfy $(3)$ either. This proves the lemma.
\end{proof}

\begin{lem}\label{L:Preimage}
The fixed point $p_\lambda$ is a hyperbolic fixed point for $f_\lambda$. If $|\lambda|\neq 0$ is sufficiently small then all the preimages 
of $p_\lambda$ in $L$ under $g_\lambda$, except $q_\lambda$ and itself, are not in $U$. And all the preimages of $q_\lambda$ under 
$g_\lambda$ are not in $U$.
\end{lem}
\begin{proof}
One readily checks that the Jacobian matrix of $f_\lambda$ at the fixed point $p_\lambda$ is
$$\left| \begin{array}{ccccc}
2t_\lambda-4 & 4 & 0 & \cdots & 0 \\
2t_\lambda-4 & 0 & 4 & \cdots & 0 \\
\vdots & \vdots & \vdots & \ddots & \vdots \\
2t_\lambda-4 & 0 & 0 & \cdots & 4 \\
2t_\lambda & 0 & 0 & \cdots & 0
\end{array} \right|$$
and that $2t_\lambda$ is one of the eigenvalues of this matrix and all the other eigenvalues 
have modulus 4. Therefore $p_\lambda$ is a hyperbolic fixed point and $L\cap U$ is the local stable manifold.

The preimages of $p_\lambda$ in $L$ under $g_\lambda$ have the form $[1:1:\cdots:1:s_\lambda]$,
with $s_\lambda$ satisfying
$$(\cdots((s_\lambda^2 + \lambda)^2 + \lambda)^2 +\cdots + \lambda)^2 + \lambda = t_\lambda$$
or equivalently
\begin{equation}
h_\lambda^k(s_\lambda) = t_\lambda,
\end{equation}
where $h_\lambda$ maps $z$ to $z^2 + \lambda$.
Note that $|\lambda|\ll\sqrt{|\lambda|}$ and $|h_\lambda^{-1}(-t_\lambda)|\sim\sqrt{|\lambda|}$,
so $|h_\lambda^{-j}(-t_\lambda)|\sim|\lambda|^{\frac{1}{2^j}}$, $1\le j\le k$. Therefore
by $(4)$, if $[1:1:\cdots:1:s_\lambda]\in U$, we must have $s_\lambda = \pm t_\lambda$.

Let us then consider the preimages of $q_\lambda$ under $g_\lambda$. There are two cases:\\
$(i)$ The preimages have the form $[1:0:\cdots:0:s_\lambda]$, with $s_\lambda$ satisfying
\begin{equation}
h_\lambda^k(s_\lambda) = -t_\lambda;
\end{equation}
$(ii)$ The preimages have the form $[z:1:\cdots:1:s_\lambda]$, with $z$ satisfying
\begin{equation}
(1-2(1-\cdots-2(1-2(1-\frac{2z^2}{(z-2)^2})^2)^2\cdots)^2)^2 = 1,
\end{equation}
and $s_\lambda$ satisfying
\begin{equation}
h_\lambda^{k-1}(\frac{s_\lambda^2 + \lambda z^2}{(z-2)^2}) = -t_\lambda.
\end{equation}

In case $(i)$, we have 
$|s_\lambda| = |h_\lambda^{-k}(-t_\lambda)|\sim|\lambda|^\frac{1}{2^k}\gg \sqrt{|\lambda|}$
by $(5)$, so $[1:0:\cdots:0:s_\lambda]\notin U$. In case $(ii)$, we have 
$\left| \frac{s_\lambda^2 + \lambda z^2}{(z-2)^2} \right| = |h_\lambda^{-(k-1)}(-t_\lambda)|\sim|\lambda|^\frac{1}{2^{k-1}}$
by $(7)$. If $z = 0$, then $|s_\lambda|\sim|\lambda|^{\frac{1}{2^k}}$, so
$[z:1:\cdots:1:s_\lambda]\notin U$. Otherwise, we can write 
$\frac{s_\lambda^2 + \lambda z^2}{(z-2)^2} = [(\frac{s_\lambda}{z})^2 + \lambda]\frac{z^2}{(z-2)^2}$.
Note that $\left| \frac{z^2}{(z-2)^2} \right| \le 1$ by $(6)$, so we have
$\left| (\frac{s_\lambda}{z})^2 + \lambda \right| \ge |h_\lambda^{-(k-1)}(-t_\lambda)|$,
$i.e.$ $\left| \frac{s_\lambda}{z} \right|\sim|\lambda|^{\frac{1}{2^k}}$.
If $|z| \ge 1$, we conclude that $[z:1:\cdots:1:s_\lambda]\notin U$. If $|z| < 1$, note
that $(6)$ only has finitely many solutions, so we can choose $\lambda$ sufficiently
small such that the smallest modulus of those solutions is greater than
$|\lambda|^{\frac{1}{2^k}}$. Then we have $|s_\lambda|\sim|\lambda|^\frac{1}{2^{k-1}}\ge \sqrt{|\lambda|}$,
and we also conclude that $[z:1:\cdots:1:s_\lambda]\notin U$. This proves the lemma.
\end{proof}

We can now prove Theorem \ref{T:Main1}. The idea of the proof is the same as in \cite{JW:Nonalgebraic}.

\begin{proof}[Proof of Theorem \ref{T:Main1}]
It suffices to prove that $M$ is nonalgebraic. And this follows if we can show that there 
is no algebraic curve $V\subset U$ which satisfies $g_\lambda(V)\subset V$. Suppose 
there is such a curve $V$. We may assume that $g_\lambda(V) = V$, because otherwise we may
replace $V$ by the algebraic curve $V' = \bigcap_{n\ge 0}g_\lambda^n(V)$. 

By Lemma \ref{L:Trapping}, the open set $U\cap L$ is mapped strictly into itself by $g_\lambda$, and is
therefore contained in the basin of attraction (for $g_\lambda|_L$) to the fixed point
$p_\lambda$. Since $g_\lambda(V)\subset V$ and $\emptyset\neq V\cap L\subset U\cap L$, we
must have $p_\lambda\in V$.

By Lemma \ref{L:Preimage}, the point $p_\lambda$ is a hyperbolic fixed point for $g_\lambda$. Since
$g_\lambda(V)\subset V$, this implies that any local irreducible branch of $V$ at
$p_\lambda$ must coincide  with either the local stable or the local unstable manifold of $g_\lambda$ 
at $p_\lambda$. But the local stable manifold is $L\cap U$, and no branch of $V$ can be 
contained in $L$. Thus the germ of $V$ at $p_\lambda$ coincides with the local unstable
manifold of $g_\lambda$ at $p_\lambda$. In particular, it is transverse to $L$ at $p_\lambda$.

If $V\cap L = \{p_\lambda\}$, then $V$ is a line by Bezout's theorem. But by Lemma \ref{L:Line}, 
there is no fixed line contained in $U$. Therefore, $V\cap L$ contains more than one point. 
By Lemma \ref{L:Preimage}, we must have $q_\lambda\in V\cap L$. Since $g_\lambda(V) = V$, the curve $V$ should contain 
at least one of the preimages of $q_\lambda$ under $g_\lambda$. But by Lemma \ref{L:Preimage}, all those 
preimages are not in $U$, contradicting $V\subset U$. This proves the nonalgebraicity of $M$, hence of $K_\lambda$.

We now show that $f_\lambda$ is chaotic on $K_\lambda$. From the proof of Corollary \ref{C:Attractor} we know
that $f_\lambda$ is topologically mixing on $K_\lambda$. So by Remark \ref{R:Mixing}, we only need to show that periodic
points for $f_\lambda$ are dense in $K_\lambda$. Let $\Omega$ be any open subset of $K_\lambda$, then
$\phi_\lambda^{-1}(\Omega)$ is an open subset of $\hat{\Pi}$, where $\phi_\lambda$ is the semiconjugation between
$\hat{\Pi}$ and $K_\lambda$ constructed in Proposition \ref{P:Semi}. Since periodic points for $f_0$ are dense in $\Pi$, 
periodic points for $\hat{f_0}$ are dense in $\hat{\Pi}$ by part $(b)$ of Lemma \ref{L:Periodic}.
Hence there exists a periodic point $\hat{x}\in \phi_\lambda^{-1}(\Omega)$. From the construction of $\phi_\lambda$ we
see that $\phi_\lambda(\hat{x})\in K_\lambda$ is a periodic point for $f_\lambda$. Therefore, periodic points for 
$f_\lambda$ are dense in $K_\lambda$. This completes the proof.
\end{proof}

\begin{rmk}
In \cite{FS:Examples}, nonalgebraic attractors for more general holomorphic maps on $\mathbf{P}^2$ are
constructed using a somewhat different method.
\end{rmk}

\section{Hyperbolic Measures}\label{S:Measure}

We study measures supported on the attractor $K_\lambda$ in this section. 

Let $\mu_0$ be the Green measure associated with $f_0|_\Pi$ and $i$ be the inclusion map that maps the hyperplane at infinity $\Pi$ into 
$\mathbf{P}^k$. Define $\mu_\lambda:={\phi_\lambda}_\star (\widehat{i_\star(\mu_0)})$, where $\phi_\lambda$ is the semiconjugation 
given in Proposition \ref{P:Semi}, and $\widehat{i_\star(\mu_0)}$ is the history measure on $\hat{\Pi}$ induced by $i_\star(\mu_0)$. 
We show that $\mu_\lambda$ satisfies Theorem \ref{T:Main2}.

\begin{proof}[Proof of Theorem \ref{T:Main2}]
Recall that the $f_0$-invariant probability measure $\mu_0$ is mixing and is the unique measure of maximal entropy $(=(k-1)\log2)$ for $f_0|_\Pi$.  
It describes the distribution of periodic points for $f_0|_\Pi$ and the smallest Lyapunov exponent of $f_0|_\Pi$ with respect to $\mu_0$ is greater 
or equal to $\frac{1}{2}\log2$ at $\mu_0-$almost every point (\cite{FS:Potential}, \cite{BrD:Exposants}, \cite{BrD:Deux}). Therefore, the measure 
$i_\star(\mu_0)$ also has the same properties. We also note from Remark \ref{R:Semi} that $\mu_\lambda={\pi_\lambda}_\star (\hat{\mu_\lambda})$, 
where $\hat{\mu_\lambda}$ is the history measure induced by  $\mu_\lambda$ and is equal to $\hat{\phi_\lambda}_\star (\widehat{i_\star(\mu_0)})$. 
Since $\hat{\phi_\lambda}$ is a conjugation between $\widehat{f_0|_\Pi}$ and $\widehat{f_\lambda|_{K_\lambda}}$, the dynamical properties of 
($\widehat{f_0|_\Pi}$, $\widehat{i_\star(\mu_0)}$) and ($\widehat{f_\lambda|_{K_\lambda}}$, $\hat{\mu_\lambda}$) are the same. Hence, 
Theorem \ref{T:Ergodic} implies $(1)$, Theorem \ref{T:Periodic} implies $(3)$, and Lemma \ref{L:Lyapunov} implies $(4)$.
From Theorem \ref{T:Entropy}, we know that the topological entropy of $f_\lambda|_{K_\lambda}$ is equal to $(k-1)\log2$. 
Hence (2) follows from Theorem \ref{T:Entropy1} and Corollary \ref{C:Entropy}. We have $supp(\mu_\lambda)=K_\lambda$ since $\phi_\lambda$ is onto. 
The fact that $\mu_\lambda$ is a probability measure invariant under $f_\lambda$ is clear from the definition and Remark \ref{R:History}.
\end{proof}

\begin{rmk}
An $f-$invariant measure $\mu$ is said to be hyperbolic if for $\mu-$almost every point $x$ the Lyapunov exponents of $f$ at $x$ are not equal to zero
(\cite{KH:Intro}). Since for $\mu_\lambda-$almost every point the smallest nonnegative Lyapunov exponent of $f_\lambda$ is positive and the Lyapunov 
exponent of $f_\lambda$ in the normal direction is negative due to the uniform contraction in that direction, the measure $\mu_\lambda$ is indeed hyperbolic.
\end{rmk}

When $k=2$, there are also other ways of defining measures supported on the attractor. We discuss two of them here.

Let $f$ be a holomorphic map on $\mathbf{P}^2$ of algebraic degree $d>1$ and $L$ be a complex line. Let $S_m=\frac{1}{m}\sum_{i=0}^{m-1}\frac{[f^i(L)]}{d^i}$
and $S$ be a cluster point of the sequence $\{S_m\}$. We have that $S$ is a positive closed $(1,1)$ current satisfying $f_{\star}S=dS$. 
Since the Green current $T$ has continuous potential and satisfies $f^{\star}T=dT$, the measure $\nu=S\wedge T$ is well-defined and is
$f$-invariant, i.e. $f_{\star}\nu=\nu$. If we start with $f=f_\lambda$ and $L=\Pi$, we get an $f_\lambda-$invariant measure $\nu_\lambda$ 
whose support is clearly contained in the attractor $K_\lambda$ since $\Pi$ is in the trapping region.

\begin{prop}\label{P:Measure}
The measure $\nu_\lambda$ is a measure of maximal entropy for $f_\lambda|_{K_\lambda}$.
\end{prop}
\begin{proof}
It is shown in \cite{deT:Selles} that the measure-theoretic entropy of $\nu_\lambda$ is at least $\log2$. On the other hand, it is shown 
in \cite{deT:Genre} that the topological entropy of $f_\lambda$ is bounded by $\log2$ outside of a neighborhood of the support of the 
Green measure $\mu$ associated with $f_\lambda$. Therefore it suffices to see that $K_\lambda$ does not intersect the support of $\mu$. 
By \cite{FS:Examples}, if $p\in supp(\mu)$ and $U$ is a neighborhood of $p$, then $\mathbf{P}^2\backslash \cup_{i=0}^{\infty}f^i(U)$ is pluripolar. 
Since $K_\lambda$ is contained in the trapping region, whose compliment is obviously not pluripolar, we must have 
$K_\lambda\cap supp(\mu)=\emptyset$.
\end{proof}

To define the other measure supported on the attractor, we need the notion of unstable manifolds. We recall the definitions here and refer
the reader to \cite{J:Hyperbolic} and \cite{S:Stability} for details. 

Let $f$ be a holomorphic map on $\mathbf{P}^2$. Let $\Gamma$ be a compact subset of $\mathbf{P}^2$ with $f(\Gamma)=\Gamma$ and
$\hat{\Gamma}$ be the history space associated with $f|_\Gamma$. For each history $\hat{p}\in \hat{\Gamma}$ and $\delta>0$ small, we define 
$local$ $unstable$ $manifolds$ by
$$W_\delta^u(\hat{p})=\{x\in \mathbf{P}^2 | \exists \hat{x}, \pi_0(\hat{x})=x, d(x_{-i},p_{-i})<\delta, \forall i\ge 0\},$$
where $d$ is the distance function on $\mathbf{P}^2$, and define $global$ $unstable$ $manifolds$ by
$$W^u(\hat{p})=\bigcup\limits_{j=0}^{\infty} f^j(W_\delta^u(\hat{f}^{-j}(\hat{p}))).$$

Now let $p$ be a hyperbolic fixed point in the attractor, and $\hat{p}$ be the history with $p_{-i}=p$ for each $i\ge 0$.
Let $D$ be the local unstable manifold associated with the history $\hat{p}$ centered at $p$. The Stable Manifold Theorem says that $D$ is an 
embedded analytic disc in $\mathbf{P}^2$. Let $\psi\ge 0$ be a test function supported in a neighborhood of $p$ and vanishing in a neighborhood 
of $\partial D$, with $\psi(p)>0$. We have the following proposition from \cite{FS:Hyperbolic}, adapted to our situation.

\begin{prop}
Let $D$ and $\psi$ be as above. Let $R_m=\frac{1}{m}\sum_{i=0}^{m-1}\frac{f_{\star}^i[\psi D]}{2^i}$, where $f=f_\lambda|_{K_\lambda}$.
Any cluster point $\tilde{R}$ of the sequence $\{R_m\}$ is a positive closed (1,1) current supported on $\overline{W^u(\hat{p})}$.
Moreover $c:=\int T\wedge \psi [D] >0$ and $R:=\frac{1}{c}\tilde{R}$ satisfies $\| R \| =1$ and $f_{\star}R=dR$.
\end{prop}

We can define then $\tau_\lambda = R\wedge T$, which is an $f_\lambda$-invariant probability measure. Note that any unstable manifold for a point in the 
attractor must be contained in the attractor by definition. Hence, the measure $\tau_\lambda$ is indeed supported on the attractor.
We show that $\tau_\lambda$ is again a measure of maximal entropy for $f_\lambda|_{K_\lambda}$. The proof is essentially the same as
in \cite{deT:Selles}, which we include for completeness.

\begin{prop}
The measure $\tau_\lambda$ is a measure of maximal entropy for $f_\lambda|_{K_\lambda}$.
\end{prop}
\begin{proof}
For simplicity, we write $f$ for $f_\lambda|_{K_\lambda}$. Denote by $B_n(x,\epsilon)$ the ball centered at $x$ with radius $\epsilon$ in the
$d_n$ metric, where $d_n(x,y)=\max_{0\le i \le n-1}\{ d(f^i(x), f^i(y)) \}$. Let $D$ and $\psi$ be as given above and define 
$\sigma_n = \frac{f^{n\star}\omega}{2^n}\wedge \frac{\psi[D]}{c}$, where $c$ is as in the previous proposition and $\omega$ is the K$\ddot{a}$hler form on 
$\mathbf{P}^2$. For any $n\ge 0$ and $\epsilon>0$ we have
$$\sigma_n(B_n(x,\epsilon)) = \frac{1}{2^n}\int_{B_n(x,\epsilon)} f^{n\star}\omega\wedge \frac{\psi[D]}{c}\le \frac{C}{2^n} v^0(f,n,\epsilon),$$
where $v^0(f,n,\epsilon)$ is the supremum of the area of $f^n(B_n(x,\epsilon)\cap D)$ over all $\epsilon$-balls and $C$ is a constant
independent of $n$. Let $\tau_n = \frac{1}{n} \sum_{i=0}^{n-1} f_{\star}^i \sigma_n$ and $\tau$ be the limit of a sub-sequence $\{\tau_{n_j}\}$.
We have the following lower-bound for the measure-theoretic entropy $h_\tau(f)$ (see \cite{deT:Selles}, \cite{BS:III})
$$h_\tau(f)\ge \limsup\limits_{n_j\rightarrow \infty} (-\frac{1}{n_j} \log(\frac{C}{2^{n_j}} v^0(f, n_j, \epsilon)))\ge \log2 - v^0(f,\epsilon),$$
where $v^0(f,\epsilon)=\limsup \frac{1}{n_j} \log^{+}(v^0(f, n_j, \epsilon))$ with $\log^{+}(a):=\max\{\log(a), 0\}$. In \cite{Y:Volume} Yomdin showed that 
$v^0(f,\epsilon)$ goes to zero as $\epsilon$ goes to zero, so that $h_\tau(f)\ge \log2$.

We next show that $\tau$ and $\tau_\lambda$ are equal.

It is easy to see that $\frac{f^{\star}\omega}{2}=\omega+dd^c u$, where $u$ is a smooth function on $\mathbf{P}^2$. Therefore
$\frac{f^{i\star}\omega}{2^i}=\omega+dd^c U_i$, where $U_i=\sum_{l=0}^{i-1} \frac{u\circ f^l}{2^l}$. Since $\frac{f^{i\star}\omega}{2^i}$
converges to the Green current $T$ (see \cite{FS:Potential}), we have $T=\omega + dd^c U$, where $U$ is a continuous function on $\mathbf{P}^2$
with $\max | U_i-U | \le \frac{C}{2^i}$. By definition
$$\tau_{n_j}=\frac{1}{n_j}\sum\limits_{i=0}^{n_j-1} f_{\star}^i(\frac{f^{n_j\star}\omega}{2^{n_j}}\wedge \frac{\psi[D]}{c})
= \frac{1}{n_j}\sum\limits_{i=0}^{n_j-1} \frac{f^{n_j-i\star}\omega}{2^{n_j-i}}\wedge \frac{f_{\star}^i(\psi[D])}{2^ic},$$
which is equal to
$$\frac{1}{n_j}\sum\limits_{i=0}^{n_j-1} (\frac{f^{n_j-i\star}\omega}{2^{n_j-i}}-T)\wedge \frac{f_{\star}^i(\psi[D])}{2^ic}
+ \frac{1}{n_j}\sum\limits_{i=0}^{n_j-1} T \wedge \frac{f_{\star}^i(\psi[D])}{2^ic}.$$
The second term is just $T\wedge \frac{R_{n_j}}{c}$ and it converges to $\tau_\lambda$ since $T$ has continuous potential. Applying the
first term to a test function $\varphi$, we get
$$\int \frac{1}{n_j}\sum\limits_{i=0}^{n_j-1} (U_{n_j-i}-U)dd^c\varphi \wedge \frac{f_{\star}^i(\psi[D])}{2^ic},$$
whose absolute value is bounded by $| \varphi |_{C^2} \frac{1}{n_j}\sum_{i=0}^{n_j-1} \frac{C}{2^{n_j-i}}$ which goes to zero as $j$ goes to
infinity. 

This shows that $\tau_\lambda=\tau$. In particular, we have $h_{\tau_\lambda}(f)\ge \log2$.

The conclusion then follows as in Proposition \ref{P:Measure}.
\end{proof}

By Theorem \ref{T:Main2}, we see that $\mu_\lambda=\nu_\lambda=\tau_\lambda$. Hence $\nu_\lambda$ and $\tau_\lambda$ have
all the properties stated in the theorem. In particular, they are supported on the whole of $K_\lambda$.  Therefore the unstable manifold $W^u(\hat{p})$ 
is dense in the attractor, a fact also noted in \cite{FS:Examples}.

\begin{rmk}
As $p$ is a hyperbolic fixed point, one of the eigenvalues of $f_\lambda^{\prime}(p)$ is strictly larger than 1 in modulus. Let $\xi$ be a
corresponding eigenvector and let $\phi: \Delta\rightarrow D$ be a holomorphic embedding with $\phi (0)= p$, where $\Delta$ is the unit disk
in $\mathbf{C}$. Since $D$ is the unstable disc associated with the history $\hat{p}$, the vector $\phi^{\prime}(0)$ is a nonzero multiple of $\xi$. Now
the sequence $\phi_n:=f_\lambda^n\circ \phi$ gives us maps from $\Delta$ into $\overline{W^u(\hat{p})}$ with $\phi_n (0)=p$ and
$|\phi_n^{\prime}(0)|$ increasing to infinity as $n$ goes to infinity. By Brody's theorem (see \cite{L:Hyperbolic}) there must be a nonconstant
entire image of $\mathbf{C}$ in $\overline{W^u(\hat{p})}$. The fact that a nontrivial attractor contains a nonconstant entire image of
$\mathbf{C}$ is valid for any holomorphic map on $\mathbf{P}^k$, $k\ge 2$ (see \cite{FW:Attractor}).
\end{rmk}

\begin{rmk}
The above two methods of constructing invariant measures supported on an attractor from currents apply to more general cases. For instance, we can
consider the following maps
$$f_{a,b}([z:w:t]) = [(z-2w)^2:z^2+at^2:t^2+bz^2],\ \ \ a, b\ small.$$
A similar computation shows that there is also a trapping region for $f_{a,b}$ containing the line at infinity. Hence the same construction as
above will give us two invariant measures of maximal entropy. But it is not clear whether a theorem similar to Theorem \ref{T:Main2}
holds true in this case. However, see \cite{Di:Attractor} for some recent results in this direction.
\end{rmk}

\end{document}